\documentclass[10pt]{amsart}%
\usepackage{enumerate}
\usepackage{amsmath}
\usepackage{amsfonts}
\usepackage{amssymb}
\usepackage{graphicx}%
\providecommand{\U}[1]{\protect\rule{.1in}{.1in}}
\newtheorem{theorem}{Theorem}
\newtheorem{corollary}[theorem]{Corollary}
\theoremstyle{remark}
\newtheorem{rem}{Remark}
\begin{document}
\title{The ideal center of the dual of a Banach lattice }
\author{Mehmet Orhon}
\address{Department of Mathematics \& Statistics, University of New Hampshire, Durham,
NH 03824}
\email{mo@unh.edu}
\dedicatory{Dedicated to Professor W.A.J. Luxemburg on the occasion of his 80th birthday}
\begin{abstract}
Let $E$ be a Banach lattice. Its ideal center $Z(E)$ is embedded naturally in
the ideal center $Z(E^{\prime})$ of its dual. The embedding may be extended to
a contractive algebra and lattice homomorphism of $Z(E)^{\prime\prime}$ into
$Z(E^{\prime})$. We show that the extension is onto $Z(E^{\prime})$ if and
only if $E$ has a topologically full center. (That is, for each $x\in E,$ the closure
of $Z(E)x$ is the closed ideal generated by $x$.) The result can be generalized
to the ideal center of the order dual of an Archimedean Riesz space and in a
modified form to the orthomorphisms on the order dual of an Archimedean Riesz space.

\end{abstract}
\date{January 31, 2010}
\maketitle

\section{Introduction}

Let $E$ be a Banach lattice and let $Z(E)$ denote its (ideal) center. In
general $Z(E)$ is a subalgebra and a sublattice of $Z(E^{\prime})$, the center
of the dual of E. The embedding may be extended to a contractive algebra and
lattice homomorphism of $Z(E)^{\prime\prime}$ into $Z(E^{\prime})$. In this
paper we will show that the extension is onto $Z(E^{\prime})$ if and only if
$E$ has a topologically full center. (That is, for each $x\in E,$ the closure
of $Z(E)x$ is the closed ideal generated by $x$.) In this case $Z(E^{\prime})$
is isomorphic, both as an algebra and as a vector lattice, to a band in
$Z(E)^{\prime\prime}$.

\indent Let $K$ be a compact Hausdorff space and let $C(K)$ denote the Banach
algebra of continuous functions on $K$ with the sup norm. Also let $L(E)$
denote the bounded linear operators on $E$ and let $m:C(K)\rightarrow L(E)$ be
a bounded unital algebra homomorphism. If each closed $C(K)$-invariant
subspace of $E$ is an ideal, we will call $m$ \textit{an ideal generating
representation} of $C(K)$ on $E$. The main result of this paper asserts that
an ideal generating representation of $C(K)$ uniquely determines the center of
both $E$ and $E^{\prime}$. Moreover, up to lattice isometry, an ideal
generating representation determines the order structure of $E$ as well.

The center of a Banach lattice has been of considerable importance in the
study of operators on Banach lattices, especially in proving dominance
theorems. Several definitions have been given in the literature of conditions
expressing that the center is large in some sense. Important examples include
Meyer's \textit{topological richness }\cite{M}, Hart's \textit{transitivity}
\cite{H} and Wickstead's \textit{topological fullness }\cite{TW1}. For our
purposes, topological fullness is the most useful of these concepts.

\indent The result stated in the first paragraph can be generalized to the
ideal center of the order dual of a Riesz space and, in a modified form, to
the orthomorphisms on the order dual of a Riesz space when the order dual
separates the points of the Riesz space.

\indent All standard terminology and standard results about vector lattices
that we use may be found in at least one of \ \cite{AB}, \cite{LZ}, \cite{S}
or \cite{Z}. Except in Remark 5, we do not make a distinction between the real
and the complex spaces. Our results hold in either case.

\section{Ideal generating representations}

Suppose a bounded unital homomorphism $m:C(K)\rightarrow L(E)$ is given. Then
we consider its \textit{Arens extension} $m^{\ast}:C(K)^{\prime\prime
}\rightarrow L(E^{\prime})$ as follows. Associated with the homomorphism $m$,
we define three bilinear maps:
\begin{align*}
C(K)\times E\rightarrow E  &  ::(a,x)\rightarrow ax:ax=m(a)(x)\\
E\times E^{\prime}\rightarrow C(K)^{\prime}  &  ::(x,x^{\prime})\rightarrow
\mu_{xx^{\prime}}:\mu_{xx^{\prime}}(a)=x^{\prime}(ax)\\
E^{\prime}\times C(K)^{\prime\prime}\rightarrow E^{\prime}  &  ::(x^{\prime
},a)\rightarrow ax^{\prime}:ax^{\prime}(x)=a(\mu_{xx^{\prime}})\;.
\end{align*}
We recall that $C(K)^{\prime\prime}$ is isomorphic to $C(S)$ with $S$
hyperstonian. Then it is routine to show that
\[
m^{\ast}:C(K)^{\prime\prime}\rightarrow L(E^{\prime}):m^{\ast}(a)(x^{\prime
})=ax^{\prime}%
\]
is a bounded unital algebra homomorphism that is $(w^{\ast},w^{\ast
}-operator)$-continuous. Also for each $a\epsilon C(K)$, $m^{\ast}(a)$ is the
adjoint in $L(E^{\prime})$ of the operator $m(a)$ in $L(E)$. For information
on Arens extensions see \cite{Ar}.

\indent Recall that $T\epsilon L(E)$ is in the center $Z(E)$ of a Banach
lattice $E$ if there is an $M>0$ such that for each $x\epsilon E_{+},$ we have
$|T(x)|\leq Mx$. (Evidently one can take $M=||T||$.) In what follows we will
use a number of known results about the center of a Banach lattice. We mention
them here briefly. An operator $T$ on $E$ is in $Z(E)$ if and only if $T$
leaves each closed ideal of $E$ invariant \cite{TW4}. If $E$ is Dedekind
complete then less is required. In this case, $T$ is in $Z(E)$ if and only if
$T$ leaves each band of $E$ invariant. This means that when $E$ is Dedekind
complete then an operator $T$ on $E$ is in $Z(E)$ if and only if $T$ commutes
with each band projection on $E.$ Each band projection is in $Z(E).$
Consequently, one has the well known result that when $E$ is Dedekind complete
then an operator $T$ is in $Z(E)$ if and only if $T$ commutes with all the
operators in $Z(E).$ Since the dual $E^{\prime}$ of a Banach lattice $E$ is
Dedekind complete, the above statements apply in particular to $Z(E^{\prime}%
)$. Finally, the definition of membership in the ideal center and the order
structure of the dual imply that an operator $T$ on $E$ is in $Z(E)$ if and
only if its adjoint $T^{\prime}$ is in $Z(E^{\prime}).$

Let $A$ be a subset of $L(E)$, we denote by $w-cl(A)$ the weak-operator
closure of $A$ in $L(E)$.

\begin{theorem}
Let $E$ be a Banach lattice, $K$ be a compact Hausdorff space and
$m:C(K)\rightarrow L(E)$ be a bounded unital algebra homomorphism. Consider
the following statements:

\begin{enumerate}
\item Each closed $C(K)$-invariant subspace of $E$ is an ideal.

\item $Z(E^{\prime})=m^{\ast}(C(K)^{\prime\prime}).$

\item $Z(E)=w-cl(m(C(K))).$
\end{enumerate}

Then $(1)\Leftrightarrow(2)\Rightarrow(3)$. Also, in this case, $m$ is a
positive contractive homomorphism.
\end{theorem}

\begin{proof}
Assume $(1)$ holds. Then, by duality, each $w^{\ast}$-closed $C(K)$-invariant
subspace of $E^{\prime}$ is a band. Therefore any operator in $Z(E^{\prime})$
leaves such a subspace invariant. But, since $C(K)$ is $w^{\ast}$-dense in
$C(K)^{\prime\prime},$ the continuity properties of $m^{\ast}$ imply that a
$w^{\ast}$-closed $C(K)$-invariant subspace of $E^{\prime}$ is in fact
$C(K)^{\prime\prime}$-invariant. We can summarize the discussion with the
statement that each $w^{\ast}$-closed $C(K)^{\prime\prime}$-invariant subspace
of \ $E^{\prime}$ is left invariant by each operator in $Z(E^{\prime}).$ Then
by a deep result due to Arenson (\cite{AAK}, 9.12, page 63), $Z(E^{\prime})$
is contained in $m^{\ast}(C(K)^{\prime\prime})$. (See Remark 1 below for a
more detailed explanation.) Since $Z(E^{\prime})$ is contained in $m^{\ast
}(C(K)^{\prime\prime})$ and since $m^{\ast}(C(K)^{\prime\prime})$ is a
commutative subalgebra of $L(E^{\prime})$, each operator in $m^{\ast
}(C(K)^{\prime\prime})$ commutes with all the operators in $Z(E^{\prime})$. It
follows from the discussion preceeding Theorem 1 that $m^{\ast}(C(K)^{\prime
\prime})$ is, in turn, a subset of $Z(E^{\prime})$. That is, $(1)\Rightarrow
(2).$

\indent Conversely, assume $(2)$ holds. Suppose $F$ is a closed $C(K)$%
-invariant subspace of $E.$ Consider its annihilator $F^{\circ}$ in
$E^{\prime}.$ $F^{\circ}$ is a $w^{\ast}$-closed $C(K)$-invariant subspace of
\ $E^{\prime}.$ Then, as in the proof of $(1)\Rightarrow(2),$ $F^{\circ}$ is a
$w^{\ast}$-closed $C(K)^{\prime\prime}$-invariant subspace of $E^{\prime}$ by
the continuity properties of $m^{\ast}.$ Now, by $(2),$ we have that
$F^{\circ}$ is a $w^{\ast}$-closed $Z(E^{\prime})$-invariant subspace of
$E^{\prime}.$ Since $E^{\prime}$ is Dedekind complete, a $Z(E^{\prime}%
)$-invariant subspace of $E^{\prime}$ is an ideal . Therefore $F^{\circ}$ is
an ideal in $E^{\prime}$. It follows from duality that $F$ is an ideal in $E$.
That is, $(2)\Rightarrow(1).$

\indent The kernel of $m^{\ast}$ is a $w^{\ast}$-closed band in $C(K)^{\prime
\prime}.$ Therefore, when $(2)$ holds, the idempotents in $C(K)^{\prime\prime
}$ are mapped onto the band projections in $Z(E^{\prime})$. Since these
generate $C(K)^{\prime\prime}$ and $Z(E^{\prime})$ respectively, $(2)$ implies
that $m^{\ast}$ is contractive and positive. In fact if $(1-e)$ denotes the
band projection in $C(K)^{\prime\prime}$ onto the kernel of $m^{\ast}$ then
$Z(E^{\prime})$ is isometrically isomorphic to $eC(K)^{\prime\prime}$ both as
an algebra and as a Banach lattice. Also, to see $(2)\Rightarrow(3),$ note
that $(2)$ and the continuity properties of $m^{\ast}$ imply that
$Z(E)=w-cl(m(C(K))).$
\end{proof}

\indent We recall that the center of a Banach lattice $E$ is isomorphic to
$C(K)$ for some compact Hausdorff space $K$ both as an algebra and a vector
lattice. Moreover the center is closed in $L(E)$ with respect to the
weak-operator topology . Also recall that a Banach lattice $E$ is said to have
a \textit{\ topologically full center} if, for each $x\epsilon E_{+},$ the
closure of $Z(E)x$ is an ideal in $E$ \cite{TW1}. This leads to:

\begin{corollary}
Let $E$ be a Banach lattice and let $i:Z(E)\rightarrow L(E)$ denote the
natural embedding of $Z(E)$ into $L(E)$. Then $E$ has a topologically full
center if and only if $i^{*}(Z(E)^{\prime\prime})=Z(E^{\prime}).$ When that
holds we may identify $Z(E^{\prime})$ with a band in $Z(E)^{\prime\prime}$
both as an algebra and as a Banach lattice.
\end{corollary}

\indent Banach lattices with a topological order unit (that is, with a
quasi-interior element) and $\sigma$-Dedekind complete Banach lattices are
examples of Banach lattices with topologically full center. Not all Banach
lattices have topologically full center. In particular a Banach lattice with
trivial center does not have a topologically full center unless the vector
lattice is one dimensional  \cite{G}, \cite{TW2}. This shows that in general
the inclusion of $i^{\ast}(Z(E)^{\prime\prime})$ in $Z(E^{\prime})$ may be
strict and that one can not expect $(3)\Rightarrow(2)$ to hold in Theorem 1.

\indent It is well known that if a Banach lattice is $\sigma$-Dedekind
complete then its center is maximal abelian in the algebra of bounded
operators on the Banach lattice, e.g., \cite{TW3}. It follows from this fact
and Corollary 2  that if the center of a Banach lattice is topologically full
then the center is maximal abelian in the algebra of bounded operators on the
Banach lattice. Also, recently Wickstead \cite{TW3} has given a direct proof
\ of the result that does not require duality arguments.

\begin{corollary}
Let $E$ be a Banach lattice with topologically full center and let $T$ be an
operator on $E$ that commutes with $Z(E)$. Then $T$ is in $Z(E).$
\end{corollary}

\begin{proof}
Consider $T^{\prime},$ the adjoint of $T$ in $L(E^{\prime})$. Then the
continuity properties of $i^{\ast}$ and of $T^{\prime}$ imply that $T^{\prime
}$ commutes with $i^{\ast}(Z(E)^{\prime\prime})=Z(E^{\prime}).$ Since
$E^{\prime}$ is Dedekind complete, $T^{\prime}\epsilon Z(E^{\prime}).$
Therefore $T\in Z(E)$ by the discussion preceeding Theorem 1.
\end{proof}

\indent Next, we show that an ideal generating representation of $C(K)$ on a
Banach lattice $E$ determines the order structure of $E$ up to lattice
isomorphism. The actual result is true more generally. In the proof below
$C_{\infty}(K)=\{f+ig:$ with $f,g\in C_{\infty}(K)_{r}\}$ where $C_{\infty
}(K)_{r}$ denotes the set of all extended continuous functions on $K$ into
$\mathbb{[-\infty},\infty],$ the two point compactification of the real
numbers. An extended continuous function on $K$ is a continuous function into
$[-\infty,\infty]$ that is finite except possibly on a nowhere dense subset of
$K$ (\cite{S}, \cite{TW3}, \cite{OhO}).

\begin{theorem}
Let $E$ be a Banach lattice and let $W$ be a cone in $E$ with the properties:

\begin{enumerate}
\item $E$ is a Banach lattice with the cone $W$ and with its original norm,

\item $E$ has the same closed ideals with respect to the cones $E_{+}$ and
$W.$
\end{enumerate}

Then there is $T\epsilon Z(E)$ with $|T|=1$ such that $T(x)=/x/_{W}$ for each
$x\epsilon E_{+}.$ (Here $/./_{W}$ denotes absolute value with respect to the
lattice structure given by the cone W.) Conversely, if $T\epsilon Z(E)$ with
$|T|=1,$ then $W=T(E_{+})$ gives a cone in $E$ such that $(1)$ and $(2)$ hold.
\end{theorem}

\begin{proof}
Note that, since $E_{+}$ and $W$ give the same closed ideals in $E,$ the ideal
center with respect to either lattice structure consists of the same
operators. In fact the order structure of the ideal center is the same in both
cases. This latter statement follows from Theorem 1 if $E$ has a topological
order unit. To see it for arbitrary $E$, apply the case of a Banach lattice
with topological order unit to the closure of each principal ideal in $E$ with
$Z(E)$ restricted to each such ideal. Therefore $Z(E)=C(K)$ for some compact
Hausdorff space $K$ and $Z(E)_{+}=C(K)_{+}$ in both lattice structures. Denote
by $I(x)$ the closed ideal generated by $x\epsilon E$ in either lattice structure.

\indent Suppose initially that $E$ has a topological order unit $u$ with
respect to $E_{+}.$ Then $/u/_{W}$ is a topological order unit with respect to
$W.$ To see this observe that $I(u)=I(/u/_{W}).$ Represent the Banach lattice
with the cone $W$ as a sublattice of $C_{\infty}(K)$ such that $/u/_{W}$ is
represented by $1$ and $C(K)$ with its lattice structure corresponds to the
ideal generated by $/u/_{W}$ (\cite{S}, III.4.5). Let $f\epsilon C(K)$ be the
function that represents $u.$ Then $|f|=1.$ Let $T\epsilon Z(E)$ be the
operator that corresponds to $\bar{f}\epsilon C(K)$. It follows from the
functional representation that $T(u)=/u/_{W}.$ For each $a\epsilon C(K)_{+},$
one has $T(au)=aT(u)=a/u/_{W}=/au/_{W}.$ Then it follows that $T(x)=/x/_{W}$
for all $x\epsilon E_{+}$ when $E$ has a topological order unit. Now repeat
the above argument for any $a\epsilon C(K),$ to show that $T(|x|)=/T(x)/_{W}$
($=/x/_{W}$) for all $x\epsilon E.$ That is, $T$ is indeed a lattice
homomorphism between the two lattice structures.

\indent In the general case let, $E=\bigcup_{x\epsilon E_{+}}I(x)$ where the
family $\{I(x):x\epsilon E_{+}\}$ is directed upwards by inclusion in the
order of $E_{+}.$ Since $I(x)=I(/x/_{W})$ for each $x\epsilon E_{+}$, by the
case when there is a topological order unit, there is $T_{x}\epsilon Z(I(x))$
such that $T_{x}(z)=/z/_{W}$ for all $z\epsilon I(x)_{+}.$ Note that $0\leq
y\leq x$ for some $x,y\epsilon E_{+}$ implies that $T_{x}$ is an extension of
\thinspace\ $T_{y}$. Therefore, without ambiguity, we define an operator $T$
of norm $1$ on $E$ by $T(z)=T_{x}(z)$ whenever $z\epsilon I(x)$ for some
$x\epsilon E_{+}$. Since $T$ preserves closed ideals, $T\epsilon Z(E).$ Also
it is clear from the definition of $T$ that $|T|=1.$ This completes the first
part of the proof.

\indent The converse on the other hand is a routine exercise. Let
$/x/_{W}=T(|x|)$ for each $x\epsilon E$ and check that $(1)$ and $(2)$ are satisfied.
\end{proof}

\begin{rem}
For the purpose of consistency and clarity we will restate Arenson's Theorem
(\cite{AAK}, 9.12) in our notation and terminology: \textit{Let }$E$\textit{
be a Banach space and let }$m:C(K)\rightarrow L(E)$\textit{ be a bounded
unital algebra homomorphism. Then an operator }$T$\textit{ on }$E^{\prime}%
$\textit{ is in }$m^{\ast}(C(K)^{\prime\prime})$\textit{ if and only if }%
$T$\textit{ leaves invariant each }$w^{\ast}$\textit{-closed }$C(K)^{\prime
\prime}$\textit{-invariant subspace of }$E^{\prime}.$ In order to see that our
restatement is faithful, note that the continuity properties of $m^{\ast}$ and
basic duality theory show that $Z(X^{\ast})$ involved in the statement of
Arenson's Theorem for Banach $C(K)$-modules is equal to $m^{\ast}%
(C(K)^{\prime\prime})$ in our setting (\cite{AAK}, 7.71, p. 46). The proof of
the theorem uses several other results in \cite{AAK} that are remarkable in
their own right. In particular the proof uses an analogue of the Factorization
Theorem of Lozanovsky on Banach function spaces \cite{Loz}, \cite{Gl} that is
proved in \cite{AAK} for Banach $C(K)$-modules. \ For an earlier version of
Arenson's Theorem in the resricted setting of Bade's Theorem see (\cite{O},
Theorem 2), for a weaker version see (\cite{HO}, Lemma 5). This latter result
states (in the same circumstances as Arenson's Theorem): \textit{If an
operator }$T$\textit{ on }$E^{\prime}$\textit{ leaves invariant each }%
$w^{\ast}$\textit{-closed }$C(K)^{\prime\prime}$\textit{-invariant subspace of
}$E^{\prime}$\textit{ and commutes with }$m^{\ast}(C(K)^{\prime\prime}%
)$\textit{ then }$T$\textit{ is in }$m^{\ast}(C(K)^{\prime\prime}).$ Its proof
is considerably simpler than that of Arenson's Theorem and it may be used to
provide an alternative proof of Theorem 1 .
\end{rem}

In the next remark we will outline the above mentioned alternative proof of
Theorem 1.

\begin{rem}
To give a different proof of Theorem 1, we only need to prove (1) implies (2)
without using Arenson's Theorem. Assume (1) in the statement of Theorem 1. Let
$a\in C(K),$ then the kernel of $m^{\ast}(a)$ is a $w^{\ast}$-closed
$C(K)^{\prime\prime}$-invariant subspace of $E^{\prime}$ and is a band.
Therefore $Z(E^{\prime})$ leaves $Ker(m^{\ast}(a))$ invariant for all $a\in
C(K).$ That is, for all $T\in Z(E^{\prime})$, $a\in C(K)$ and $x^{\prime}\in
E^{\prime}$,
\[
m^{\ast}(a)(x^{\prime})=0\text{ implies }m^{\ast}(a)(T(x^{\prime}))=0.\text{ }%
\]
Then by a result of Evans \cite{E}, each operator in $Z(E^{\prime})$ commutes
with all the operators in $m^{\ast}(C(K)).$ (See also \cite{AAK}, 9.5.) Since
$E^{\prime}$ is Dedekind complete, we have $m^{\ast}(C(K))\subset Z(E^{\prime
}).$ By the Hahn-Banach Theorem, the set of non-negative elements of the unit
ball of $C(K)$ are $w^{\ast}$-dense in the set of non-negative elements of the
unit ball of $C(K)^{\prime\prime}.$ Now, the continuity properties of the
homomorphism $m^{\ast}$ and the definition of membership in the ideal center
imply that $m^{\ast}(C(K)^{\prime\prime})\subset Z(E^{\prime}).$ This means
that $Z(E^{\prime})$ leaves invariant each $w^{\ast}$-closed $C(K)^{\prime
\prime}$-invariant subspace of $E^{\prime}$ and commutes with $m^{\ast
}(C(K)^{\prime\prime}).$ Then Lemma 5 \cite{HO} implies that $Z(E^{\prime
})\subset m^{\ast}(C(K)^{\prime\prime}).$ This completes the proof.
\end{rem}

\begin{rem}
It is natural to ask if the converse of Corollary 3 is true. A.W. Wickstead
has given an example of a Banach lattice that does not have a topologically
full center, but whose center is maximal abelian \cite{TW3}. Wickstead's paper
contains several interesting problems related to Banach lattices with maximal
abelian center or with topologically full center.
\end{rem}

\begin{rem}
When $E$ is a Banach lattice with topological order unit \ a much stronger
result than Corollary 3 is proved in (\cite{AAK}, 9.11, (1)): \textit{Let }%
$T$\textit{ be an operator on }$E^{\prime}$\textit{ that commutes with
}$i^{\ast}(Z(E)),$\textit{ then }$T$\textit{ is in }$Z(E^{\prime}).$ It is
further stated in \cite{AAK} that the result remains true if $E$ is a Banach
lattice with topologically full center. However this latter statement is still
an open problem. For a more general version of this problem, the interested
reader should see (\cite{HO}, 3., p. 357). It is worth noting that the other
two questions raised in (\cite{HO}, 3.) have been answered affirmatively by
Arenson's Theorem.
\end{rem}

\begin{rem}
Let $E$ be a Riesz space with point separating order dual $E^{\sim}.$ We will
say that $E$ has a topologically full center, if for each $x\epsilon E_{+},$
the weak closure of $Z(E)x$ is an ideal. Let $L_{b}(E)$ denote the order
bounded operators on $E.$ Let $i:Z(E)\rightarrow L_{b}(E)$ be the natural
embedding and let $i^{\sim}:Z(E)^{\prime\prime}\rightarrow L_{b}(E^{\sim})$ be
its Arens extension. \textit{Then }$E$\textit{ has topologically full center
if and only if \thinspace\ }$i^{\sim}(Z(E)^{\prime\prime})=Z(E^{\sim}).$ Let
$Orth(E)$ denote the orthomorphisms on $E.$ (That is, the band preserving
operators in $L_{b}(E).$) Let $E_{n}^{\sim}$ denote the order continuous
linear functionals on $E.$ Let $\gamma:Orth(E)\rightarrow L_{b}(E)$ be the
natural embedding and let its Arens extension be denoted by $\gamma^{\sim
}:(Orth(E)^{\sim})_{n}^{\sim}\rightarrow L_{b}(E^{\sim}).$ \textit{Then }%
$E$\textit{ has topologically full center if and only if }$\gamma^{\sim
}((Orth(E)^{\sim})_{n}^{\sim})$\textit{ is an ideal in }$Orth(E^{\sim}).$ In
the case of orthomorphisms we do not need the bidual since $(Orth(E)^{\sim
})_{n}^{\sim}$ is the bidual of $Orth(E)$ \cite{HBP}. The motivation for the
results stated in this remark is the work of Huijsmans and de Pagter
\cite{HBP} on the Arens product on the bidual of an f-algebra. The above
statements are proved using the results developed in \cite{HBP} and elementary
duality theory. \newline
\end{rem}

%\begin{acknowledgement}
I wish to thank A.W. Wickstead for his valuable comments and the referee for
suggesting the inclusion of the alternative proof of Theorem 1 and the remarks
on the open questions.
%\end{acknowledgement}

\end{document}